\documentclass[10pt]{amsart}
\usepackage{amsmath,amssymb,amsthm}
\usepackage[margin=30mm]{geometry}
\usepackage[colorlinks=true,linkcolor=red,citecolor=blue,urlcolor=blue]{hyperref}
\numberwithin{equation}{section}
\newtheorem{theorem}{Theorem}[section]
\newtheorem{lemma}[theorem]{Lemma}
\newtheorem{proposition}[theorem]{Proposition}
\newtheorem{corollary}[theorem]{Corollary}
\theoremstyle{remark}
\newtheorem{remark}[theorem]{Remark}
\newcommand{\R}{\mathbb R}
\def \L{{\mathcal L}}
\def \H{{\mathcal H}}
\def \C{{\mathcal C}}
\def \N{{\mathbb N}}
\newcommand{\dd}{\,dr}

\title[Sharp low-dimensional Hydrogen uncertainty principle]{Sharp low-dimensional Constants in the second-order hydrogen Heisenberg uncertainty principle via Jacobi chains}
\author{Xia Huang}
\address{School of Mathematical Sciences, Key Laboratory of MEA (Ministry of Education) and Shanghai Key Laboratory of PMMP, East China Normal
University, Shanghai, China}
\email{xhuang@cpde.ecnu.edu.cn}
\author{Dong Ye}
\address{School of Mathematical Sciences, Key Laboratory of MEA (Ministry of Education) and Shanghai Key Laboratory of PMMP, East China Normal
University, Shanghai, China}
\email{dye@math.ecnu.edu.cn}

\keywords{Hydrogen uncertainty principle, sharp constant, Laguerre polynomials, Jacobi operator, continued fraction, Coulomb operator}
\begin{document}
\raggedbottom
\begin{abstract}
We establish the sharp constants $\H_N$ of the second-order Hydrogen Heisenberg uncertainty principle for $N=2,3$. Starting from the first angular mode reduction, we rewrite the fourth-order radial quotient as two coupled radial Coulomb channels. Then we turn the compatibility constraint into an explicit Jacobi chain by expanding the $B$-channel in generalized Laguerre polynomials and passing to parity-tail variables. Its tail energies satisfy scalar Riccati recursions, equivalently continued fractions. Furthermore, we characterize $2\H_2$ and $2\H_3$ as the unique roots of two explicit scalar equations, and we show that free-tail and Dirichlet finite sections provide monotone lower and upper approximations converging to the sharp constants. This yields rigorous bounds and the values 
$$\H_2\approx0.84478056,\quad \H_3\approx1.82348779.$$
We also obtain uniqueness of the first-mode radial extremizer up to amplitude, rotation, and dilation.
\end{abstract}
\maketitle

\section{Introduction}
For $N\ge2$, let
\begin{equation}\label{eq:1.1}
\H_N:=\inf_{u\in W^{2,2}(\R^N),\ u\ne0}
\frac{\|\Delta u\|_2\|\nabla u\|_2}{\displaystyle\int_{\R^N}\frac{|\nabla u|^2}{|x|}\,dx}.
\end{equation}
The sharp second-order hydrogen Heisenberg uncertainty principle was introduced and studied systematically by Cazacu, Flynn and Lam \cite{CFL}. They established the value $\frac{N+1}{2}$ for $N\geq 5$ and conjectured that the same value should remain valid in the remaining dimensions. Chen and Tang \cite{CT} exhibited symmetry breaking in $N=2,3$. The authors \cite{HY} settled the conjectured case $N=4$ and improved the low-dimensional bounds. In particular, we gave
\begin{equation}\label{eq:1.2}
\frac{3+6\sqrt2}{14}\le \H_2\le\frac{\sqrt3}{2},\quad
\frac74\le \H_3\le\frac{2\sqrt{21}}5,
\end{equation}
while the exact values of $\H_2$ and $\H_3$ were left open.

The aim of this  paper is to convert the two remaining low-dimensional variational problems into explicit one-dimensional Jacobi spectral problems. The starting point is the linearization principle used in \cite{HY}: after a scaling optimization,
\begin{equation}\label{eq:1.3}
\H_N^+:=\inf_{u\ne0}\frac{\|\Delta u\|_2^2+\|\nabla u\|_2^2}{\displaystyle\int_{\R^N}|x|^{-1}|\nabla u|^2\,dx}=2\H_N.
\end{equation}
For $N=2,3$, spherical harmonics reduce the unresolved part of \eqref{eq:1.3} to the first angular mode and to a fourth-order one-dimensional quotient. Our first observation is that this quotient admits a hidden two-channel factorization. The channels are radial Coulomb operators with spectra given by generalized Laguerre polynomials. The channels are not independent: they are coupled by the fact that the original vector field is a gradient. Nevertheless, after passing to Laguerre coefficients and then to parity-tail variables, the compatibility constraint becomes nearest-neighbor.

The resulting of Jacobi chains are explicit. Their effective tail energies satisfy scalar Riccati recursions, hence continued fractions. This gives exact spectral characterizations of both constants, as well as canonical lower and upper finite-section schemes.
\begin{theorem}\label{thm:main}

For $N = 2,3$, we have $2\H_N=\alpha_N$, where $\alpha_2\in(1,\sqrt3)$ and $\alpha_3\in(2,4)$ are respectively the unique solution of
\begin{equation}\label{eq:1.4}R_0(\alpha)=\frac{\alpha-1}{2}\quad \mbox{and}\quad \overline R_0(\alpha)=\frac{\alpha-2}{8},\end{equation}
here $R_0$, $\overline R_0$ are given by the explicit continued fraction in \eqref{eq:7.3}, \eqref{eq:7.4} respectively.

Moreover, the free-tail finite sections give an increasing sequence of lower approximations, and the Dirichlet finite sections give a decreasing sequence of upper approximations. Both sequences converge to the corresponding $\alpha_N$. In particular,
\begin{equation}\label{eq:1.6} \H_2\approx0.84478056,\quad \H_3\approx1.82348779.\end{equation}
\end{theorem}

Notice that Theorem~\ref{thm:main} gives an exact spectral characterization of the best constants $\H_2$ and $\H_3$, not the decimal approximation. More precisely, the constants are determined as unique spectral thresholds of explicit Jacobi pencils, and the finite sections give an algorithm for computing them to arbitrary precision. In the following, we use different notations for the sharp values, and we have
$$\alpha_N := \C_N = 2\H_N.$$

The paper is organized as follows. Section~2 recalls the low-dimensional radial reduction. Section~3 treats the two-dimensional case, with Coulomb splitting and Jacobi chain, the continued-fraction with finite-section approximations, and the extremizer representation. Section~4 gives the three-dimensional analogue in the energy completion. The appendix contains the elementary algebra behind the edge coefficients.

\section{The first-mode reduction}
We use the spherical harmonic expansion
\[
u(x)=\sum_{k=0}^\infty r^k v_k(r)\phi_k(\sigma),\quad \sigma=\frac{x}{|x|},
\]
where $-\Delta_{\mathbb S^{N-1}}\phi_k=k(N+k-2)\phi_k$. The mode estimates in \cite{HY} imply that in dimensions $N=2,3$ the only unresolved angular contribution is $k=1$. More precisely, the radial mode has quotient $N+1$, and all modes $k \ge 2$ have quotient at least $N+1$. By the symmetry-breaking examples of Chen-Tang, the first-mode quotient $\C_N$ is strictly less than $N+1$, hence
\begin{equation}\label{eq:2.1}
2\H_N= \C_N:=\inf_{v\in C_c^2(0,\infty),\ v\ne0}\frac{E_N(v)}{F_N(v)},\quad \mbox{for }\; N=2,3,
\end{equation}
where
\begin{equation}\label{eq:2.2}
E_N(v)=\int_0^\infty r^{N+1}v''^2(r)\dd+(N+1)\int_0^\infty r^{N-1}v'^2(r)\dd+\int_0^\infty r^{N+1}v'^2(r)\dd,
\end{equation}
\begin{equation}\label{eq:2.3}
F_N(v)=\int_0^\infty r^Nv'^2(r)\dd+\int_0^\infty r^{N-2}v^2(r)\dd.
\end{equation}
These are exactly the functionals in \cite[(3.9)]{HY}. We shall work first on $C_c^2(0,\infty)$ and then pass to the energy completion. Define
\begin{equation}\label{eq:2.4}G_{N,\alpha}(v):=E_N(v)-\alpha F_N(v).\end{equation}
Thus $\C_N$ is the nonnegativity threshold for $G_{N,\alpha}$.

\section{Dimension two}
We determine $\C_2 = 2\H_2$ in four steps. First, we rewrite the fourth-order radial quotient \eqref{eq:2.1}-–\eqref{eq:2.3} in terms of two Coulomb channels $A$ and $B$ (Section 3.1). Second, we expand the $B$-channel in generalized Laguerre functions, express the gradient constraint $A= T_2B$ as a parity-tail relation, and reduce the problem to an explicit Jacobi chain (Section 3.2). Third, we solve the chain by a backward Riccati recursion, equivalently by a continued fraction; free-tail and Dirichlet finite sections give monotone lower and upper approximations to $\C_2$ (Section 3.3). Finally, we reconstruct the extremizer from the minimizing Jacobi sequence (Section 3.4).

\subsection{Two Coulomb channels}
Fix $N=2$, and set
\begin{equation}\label{eq:3.1}B(r):=\frac r2v'(r),\quad A(r):=v(r)+B(r).\end{equation}
Then
\begin{equation}\label{eq:3.2}v=A-B,\quad rv'=2B,\quad rA'=rB'+2B.\end{equation}
\begin{lemma}\label{lem:3.1}
For $v\in C_c^2(0,\infty)$ and $A,B$ defined by \eqref{eq:3.1}, then there holds
\begin{equation}\label{eq:3.3} \frac{F_2(v)}2 =\int_0^\infty(A^2+B^2)\dd,\end{equation}
and
\begin{equation}\label{eq:3.4}\frac{E_2(v)}2 =\int_0^\infty r(A'^2+A^2)\dd+\int_0^\infty\Big(rB'^2+\frac4rB^2+rB^2\Big)\dd.\end{equation}
\end{lemma}
\begin{proof}
The identity \eqref{eq:3.3} follows from
\[2(A^2+B^2)=2v^2+2rvv'+r^2v'^2 = r^2v'^2 + v^2 + \big(rv^2\big)'.\]
Since
\[A'=\frac32v'+\frac r2v'',\quad B'=\frac12v'+\frac r2v'',\]
expansion of the right-hand side of \eqref{eq:3.4}, followed by integration by parts in the mixed terms, gives the identity \eqref{eq:3.4}.
\end{proof}
Define two quadratic forms as follows:
\begin{equation}\label{eq:3.5}Q_0[A]:=\int_0^\infty r(A'^2+A^2)\dd \quad\mbox{and} \quad Q_2[B]:=\int_0^\infty\Big(rB'^2+\frac4rB^2+rB^2\Big)\dd.\end{equation}
The relationship between $A$ and $B$ is given as follows.
\begin{lemma}\label{lem:3.2}
If $A-B$ vanishes at infinity and \eqref{eq:3.2} holds, then
\begin{equation}\label{eq:3.7} A(r)=T_2B(r):=B(r)-2\int_r^\infty\frac{B(s)}s\,ds := B(r) + R_2B(r).\end{equation}
Moreover $T_2$ extends continuously from the form domain of $ Q_2$ into the form domain of $Q_0+\|\cdot\|_{L^2(dr)}^2$.
\end{lemma}
\begin{proof}
The formula \eqref{eq:3.7} follows from $(A-B)'=2r^{-1}B$. Let now $Q_2(B) < \infty$.
The weighted Hardy estimate provides
\begin{equation}\label{eq:3.8}
\int_0^\infty r^a|R_2B(r)|^2\dd\le\frac {16}{(a+1)^2}\int_0^\infty r^a|B(r)|^2\dd, \quad \mbox{for any $a>-1$} .
\end{equation}
With $a=0,1$, we control $\|R_2B\|_{L^2((1+r)dr)}$. Since $A'=B'+2r^{-1}B$, $\|A'\|_{L^2(rdr)}^2 \le CQ_2[B]$.
\end{proof}
The two channel operators associated to $Q_0$ and $Q_2$ are respectively
\[\L_0 A=-(rA')'+rA,\quad \L_2 B=-(rB')'+\Big(\frac4r+r\Big)B.\]
Let $L_n^{(\beta)}$ denote the generalized Laguerre polynomial for $\beta > -1$, more details see \cite{Szego}
$$L_n^{(\beta)}(x)
=\sum_{k=0}^{n}(-1)^k\binom{n+\beta}{n-k}\frac{x^k}{k!}
=\frac{(\beta+1)_n}{n!}\,{}_1F_1(-n;\beta+1;x).$$
It is well known that
$$\int_0^\infty x^\beta e^{-x}L_n^{(\beta)}(x)L_m^{(\beta)}(x)\,dx
=\frac{\Gamma(n+\beta+1)}{n!}\delta_{nm}.$$
Denote $L_n^{(0)}=L_n$ and define
\begin{equation}\label{eq:3.9}\phi_n(r)=e^{-r}L_n(2r),\quad \psi_m(r)=r^2e^{-r}L_m^{(4)}(2r).\end{equation}
The Laguerre equation
\begin{align*}
x L''_n(x) + (\beta+1-x)L'_n(x) + n L_n(x)=0
\end{align*}
 gives
\begin{equation}\label{eq:3.10} \L_0\phi_n =(2n+1)\phi_n,\quad \L_2\psi_m =(2m+5)\psi_m,\end{equation}
\begin{equation}\label{eq:3.11} \int_0^\infty\phi_n\phi_k\dd=\frac{\delta_{nk}}{2},\quad\int_0^\infty\psi_m\psi_\ell \dd=h_m\delta_{m\ell},\end{equation}
where $h_m=\frac{(m+1)(m+2)(m+3)(m+4)}{32}.$
\begin{lemma}\label{lem:3.3}
For every $m\ge0$, there holds
\begin{equation}\label{eq:3.13}
T_2\psi_m =-\frac{2m+5}2\sum_{\substack{0\le n\le m\\n\equiv m\ (\mathrm{mod}\ 2)}}\phi_n +\frac{(m+1)(m+2)}4\phi_{m+2}.
\end{equation}
\end{lemma}
\begin{proof} The proof is a direct comparison of the generating functions for $L_n$ and $L_m^{(4)}$ after applying $T_2$. More precisely,
the Laguerre generating functions for $L_n$ and $L_m^{(4)}$ are
\begin{equation*}\sum_{n\ge0}L_n(2r)z^n=(1-z)^{-1}\exp\Big({-\frac{2rz}{1-z}}\Big),\end{equation*}
and
\begin{equation*}\sum_{m\ge0}L_m^{(4)}(2r)z^m=(1-z)^{-5}\exp\Big({-\frac{2rz}{1-z}}\Big).\end{equation*}
Set $a(z)= \frac{1+z}{1-z}$. Multiplying by the exponential factors in \eqref{eq:3.9} gives
\[\sum_{m\ge0}\psi_m(r)z^m=r^2(1-z)^{-5}e^{-a(z)r}.\]
Since
\[\int_r^\infty se^{-as}\,ds=e^{-ar}\Big(\frac ra+\frac1{a^2}\Big),\]
application of $T_2$ gives
\begin{equation}\label{eq:A.3}\sum_{m\ge0} T_2\psi_m z^m =(1-z)^{-5}e^{-a(z)r}\Big[r^2-\frac{2r}{a(z)}-\frac2{a(z)^2}\Big].\end{equation}
On the other hand,
\[\sum_{n\ge0}\phi_n(r)z^n=(1-z)^{-1}e^{-a(z)r}.\]
Summing the right-hand side of \eqref{eq:3.13} over $m$ and evaluating the even and odd parity tails gives exactly \eqref{eq:A.3}. Equality of coefficients of $z^m$ proves \eqref{eq:3.13}.\end{proof}
Expand
\begin{equation*}
B=\sum_{m\ge0}b_m\psi_m\end{equation*}
and introduce parity-tail variables
\begin{equation}\label{eq:3.15}s_n:=\sum_{\substack{m\ge n\\m\equiv n\ (\mathrm{mod}\ 2)}}\frac{2m+5}2b_m.\end{equation}
Then
\begin{equation}\label{eq:3.16}b_n=\frac2{2n+5}(s_n-s_{n+2}).\end{equation}
If $A=\sum a_n\phi_n$, Lemma~\ref{lem:3.3} yields
\begin{equation}\label{eq:3.17}a_n=\kappa_ns_{n-2}-(1+\kappa_n)s_n,\quad\kappa_n=\frac{n(n-1)}{2(2n+1)},\end{equation}
with $s_{-2}=s_{-1}=0$.
\begin{proposition}\label{prop:3.4}
Finite linear combinations of $\{\psi_m\}$ form a core for the $B$-channel. Under $T_2$ their images converge in the $A$-channel form norm. Consequently the finite Laguerre coefficient reductions are dense in the original $(E_2+F_2)^{1/2}$ energy space.
\end{proposition}
\begin{proof}
Finite Laguerre sums are dense in the $Q_2$-form domain. As $T_2$ is bounded from the $Q_2$-form domain into the $Q_0$-form domain, their images are dense in the graph norm. Lemma~\ref{lem:3.1} identifies this graph norm with the original one-dimensional energy.
\end{proof}

\subsection{The Jacobi chain}
By \eqref{eq:3.15}, we can decouple the parity sectors. Put
\begin{equation*}x_j=s_{2j},~~y_j=s_{2j+1},\quad j\ge0,\end{equation*}
and define
\begin{equation*}\mathcal E_j(x,y)= p_jx^2-2q_j xy+r_j y^2,\quad \mathcal E'_j(x,y)= p'_jx^2-2q'_j xy+r'_j y^2
\end{equation*}
where
\begin{equation}\label{eq:4.3}
p_j =\frac{(2j+1)(2j+2)(4j^2+10j+7)}{4(4j+5)^2}, \quad q_j =\frac{(2j+1)(2j+2)(2j+3)(2j+4)}{4(4j+5)^2},
\end{equation}
\begin{equation}\label{eq:4.4}
p'_j =\frac{(2j+2)(2j+3)(4j^2+14j+13)}{4(4j+7)^2}, \quad q'_j =\frac{(2j+2)(2j+3)(2j+4)(2j+5)}{4(4j+7)^2},
\end{equation}
and
 \begin{equation}\label{eq:4.5}r_j =\frac{(2j+3)(2j+4)(4j^2+10j+7)}{4(4j+5)^2}, \quad r'_j=\frac{(2j+4)(2j+5)(4j^2+14j+13)}{4(4j+7)^2}.
\end{equation}
A direct simplification gives
\begin{equation}\label{eq:4.6}p_j r_j - q_j^2=\frac{(2j+1)(2j+2)(2j+3)(2j+4)}{16(4j+5)^2}>0,
\end{equation}
and
\begin{align}\label{eq:4-6}
p'_j r'_j - (q'_j)^2=\frac{(2j+2)(2j+3)(2j+4)(2j+5)}{16(4j+7)^2}>0.
\end{align}
\begin{theorem}\label{thm:4.1}
For $\alpha \in\mathbb{R}$, we have $G_{2,\alpha}(v)=G_{2,\alpha,\rm odd}(v)+ G_{2,\alpha,\rm even}(v)$, with
\begin{equation}\label{eq:4.7}\frac{G_{2,\alpha,\rm even}(v)}2=\frac{1-\alpha}2x_0^2+\sum_{j=0}^\infty(4j+5-\alpha)\mathcal E_j (x_j,x_{j+1}),\end{equation}
and
\begin{align}\label{eq:4-7}
\frac{G_{2,\alpha,\rm odd}(v)}2=\frac{3-\alpha}2y_0^2+\sum_{j=0}^\infty(4j+7-\alpha)\mathcal E'_j (y_j,y_{j+1}).
\end{align}
Consequently
\begin{equation}\label{eq:4.8} \C_2=\inf_{x\ne0}\frac{\frac12 x_0^2+\sum_{j\ge0}(4j+5)\mathcal E_j (x_j,x_{j+1})}{\frac12x_0^2+\sum_{j\ge0}\mathcal E_j (x_j,x_{j+1})}.\end{equation}
\end{theorem}
\begin{proof}
By \eqref{eq:3.10}--\eqref{eq:3.11},
\[\frac{G_{2,\alpha}(v)}2=\frac12\sum_{n\ge0}(2n+1-\alpha)a_n^2+\sum_{n\ge0}(2n+5-\alpha)h_n b_n^2.\]
As mentioned, we can separate the study of odd and even sectors thanks to \eqref{eq:3.15}. Using \eqref{eq:3.16} and \eqref{eq:3.17}, we can expand $G_{2,\alpha}$ in function of $s_j$. For even $n=2j$, $a_0=-x_0$, while $b_{2j}$ and $a_{2j+2}$ share the factor $4j+5-\alpha$. Substituting \eqref{eq:3.16} and \eqref{eq:3.17} gives \eqref{eq:4.7}; since the pair $(b_n,a_{n+2})$ contributes
\begin{multline*}
(2n+5-\alpha)\Big\{\frac{(n+1)(n+2)(n+3)(n+4)}{8(2n+5)^2}(x_j-x_{j+1})^2\\
+\frac12\Big[\frac{(n+1)(n+2)}{2(2n+5)}x_j-\Big(1+\frac{(n+1)(n+2)}{2(2n+5)}\Big)x_{j+1}\Big]^2\Big\}.
\end{multline*}
Hence we get \eqref{eq:4.6}. Similarly, considering the odd sector with $y_j$, we have \eqref{eq:4-6} and \eqref{eq:4-7}.

By \cite{HY}, we need only to consider $0<\alpha<3$, since $2 \H_2\leq \sqrt 3<3$. Fix now $0<\alpha<3$, there holds $G_{2,\alpha,\rm odd}\geq 0$, so the critical value of $\alpha$ such that $G_{2,\alpha}\geq 0$ coincides with the largest value $\alpha$ satisfying $G_{2,\alpha,\rm even}\geq 0$. Moreover, Proposition~\ref{prop:3.4} justifies passage from finite to general form-domain sequences, so we need only to consider the crucial $\alpha \in (0,3)$ such that the right hand side of \eqref{eq:4.7} is nonnegattive for all sequences $x_j$. Isolating the coefficients of $\alpha$ yields \eqref{eq:4.8}.
\end{proof}

\subsection{Continued fractions and monotone finite sections}
Recall that $N=2$, write
\[\lambda_j =4j+5,\quad \beta_2(\alpha)=\frac{\alpha-1}2.\]
Let $p_j,q_j,r_j$ denote the edge coefficients from \eqref{eq:4.3}--\eqref{eq:4.5}. For $\alpha<\lambda_0 = 5$, define the effective tail energy
\begin{equation}\label{eq:7.1}
R_j(\alpha):=\inf_{\substack{x_j=1\\ x\text{ in the chain form domain}}}\sum_{k=j}^\infty(\lambda_k-\alpha)\mathcal E_k(x_k,x_{k+1}).
\end{equation}
The strict positivity of $\mathcal E_k$ implies $R_j(\alpha)>0$.
\begin{proposition}\label{prop:7.1}
For $N=2$, there holds the iterative formula
\begin{equation}\label{eq:7.2}
R_j(\alpha)=(\lambda_j-\alpha)p_j-\frac{(\lambda_j-\alpha)^2q_j^2}{(\lambda_j-\alpha)r_j+R_{j+1}(\alpha)}.
\end{equation}
\end{proposition}
\begin{proof}
Fix $x_j=1$, isolate the first edge and write $t=x_{j+1}$. The quantity to minimize is
\[(\lambda_j-\alpha)\bigl(p_j-2q_jt+r_jt^2\bigr)+R_{j+1}(\alpha)t^2.\]
It is strictly convex in $t$, minimizing it gives \eqref{eq:7.2}.
\end{proof}
Iteration of \eqref{eq:7.2} gives
\begin{equation}\label{eq:7.3}
R_0(\alpha)=c_0p_0 -\cfrac{c_0^2q_0^2}{c_0r_0 +c_1p_1 -\cfrac{c_1^2q_1^2}{c_1r_1 +c_2p_2-\cdots}},\quad c_j=4j+5-\alpha.
\end{equation}
\begin{theorem}\label{thm:7.2}
The constant $\C_2=2\H_2$ is characterized by
\begin{equation}\label{eq:7.5}R_0(\C_2)=\beta_2(\C_2).\end{equation}
The equation has at most one solution in the interval $(-\infty, 5)$. Hence the solution in Theorem~\ref{thm:main} is unique.
\end{theorem}
\begin{proof}
Normalize $x_0=1$. The Jacobi representations \eqref{eq:4.7} give
\[\inf_{x_0=1}\frac{G_{2,\alpha,\rm even}}2=-\beta_2(\alpha)+R_0(\alpha),\]
A negative direction cannot occur with $x_0=0$ as long as all edge factors are positive. Thus the nonnegativity threshold is exactly \eqref{eq:7.5}.

For every fixed admissible tail, the energy in \eqref{eq:7.1} is nonincreasing in $\alpha$. Hence $R_0(\alpha)$ is nonincreasing, while $\beta_2(\alpha)$ is strictly increasing. Therefore the threshold equation has at most one solution. Existence follows from the variational threshold $\C_2$ and the equivalence of the chain and continuous form domains established above.
\end{proof}
We next make the approximation statement precise. Fix $J\ge0$. The \emph{free-tail} approximation discards all edges with index $k\ge J$, while the \emph{Dirichlet} approximation restricts to $x_k=0$ for every $k\geq J+1$. Equivalently, in the backward Riccati recursion use
\begin{equation*}R_J^{-,J}=0,\quad R_J^{+,J}=(\lambda_J-\alpha)p_J.\end{equation*}
Let $\alpha_J^\pm$ be the roots of
\[R_0^{\pm,J}(\alpha)=\beta_2(\alpha).\]
\begin{proposition}\label{prop:7.3}
For any $J$,
\begin{equation*}\alpha_{J}^-\le\alpha_{J+1}^-\le \C_2\le\alpha_{J+1}^+\le\alpha_{J}^+,\quad \lim_{J\to\infty}\alpha_{J}^- = \lim_{J\to\infty}\alpha_{J}^+ = \C_2.\end{equation*}
\end{proposition}
\begin{proof}
The free-tail problem is obtained by deleting nonnegative edge terms. Adding one further positive edge can only increase its normalized effective energy, hence the lower approximations increase and remain below the full tail energy. Conversely, Dirichlet truncation is a restriction to a finite-dimensional subspace; enlarging the subspace decreases the restricted minimum, so the upper approximations decrease and remain above the full value.

For the lower sequence, the increasing finite positive forms converge in the sense of monotone closed-form convergence to the infinite Jacobi form. Equivalently, one may take finite-tail minimizers normalized by $x_0=1$, extract a coordinatewise weakly convergent subsequence using the coercivity of the positive edge forms, and apply lower semicontinuity. For the upper sequence, finite-support sequences are dense by Proposition~\ref{prop:3.4}; hence the Rayleigh--Ritz principle gives convergence of the Dirichlet finite sections. The uniqueness in Theorem~\ref{thm:7.2} identifies limit with $\C_2$.
\end{proof}
\begin{remark}\label{rem:7.4}
All coefficients of the finite Jacobi matrices are rational when $\alpha$ is rational. Consequently positive definiteness or indefiniteness can be checked by the exact $LDL^T$ recursion
\[d_0=D_0,\quad d_j=D_j-\frac{E_{j-1}^2}{d_{j-1}}.\]
Using the exact rational arithmetic (for example with Python), we obtain the rigorous bounds
\[0.844780555\le \H_2<0.844780565\]
The exact value is given by the continued-fraction characterization.
\end{remark}

\subsection{Extremizers and reconstruction}
At the critical value $\alpha=\C_2$, normalize $x_0=1$. The tail minimization is strictly convex. Hence the minimizing Jacobi sequence is determined recursively by
\begin{equation}\label{eq:8.1}
x_{j+1}=\frac{(\lambda_j-\C_2)q_j}{(\lambda_j-\C_2)r_j +R_{j+1}(\C_2)}x_j.
\end{equation}
Thus the reduced extremizer is unique up to a scalar multiple. The $B$-channel coefficients are
\begin{equation*}b_{2j} =\frac{2(x_j-x_{j+1})}{4j+5},\quad b_{2j+1}=0,\end{equation*}
and
\begin{equation*}B^*(r)=\sum_{j\ge0}b_{2j} r^2e^{-r}L_{2j}^{(4)}(2r),\quad v^*(r)=-2\int_r^\infty\frac{B^*(s)}s\,ds.\end{equation*}
The first spherical harmonic then reconstructs the linearized extremizers as
\begin{equation}\label{eq:8.6}u^*(x)=(a\cdot x)v^*(|x|),\quad a\in\R^2.\end{equation}
Amplitude, rotation, and the dilation arising when passing from the linearized inequality back to the product inequality give the full natural extremal family.
\begin{corollary}
The first-mode radial profile is unique up to multiplication by a nonzero scalar. Consequently the extremal family is unique modulo amplitude, rotation and dilation.
\end{corollary}
\begin{proof}
The denominator in \eqref{eq:8.1} is strictly positive and each one-step minimization is strictly convex. Therefore $x_0$ uniquely determines the entire minimizing sequence. The Laguerre reconstruction and the first spherical harmonic then give \eqref{eq:8.6}. The remaining symmetries are those of the original product inequality.
\end{proof}

\section{Dimension three}
The approach is very similar to the two-dimensional case, so we will go through quickly. For $N=3$, set
\begin{equation}\label{eq:5.1} A:=v +\frac r3v',\;\; B:=\frac{\sqrt2}3rv', \quad \mbox{hence }\; \sqrt2 rA'=rB'+3B.\end{equation}
Then direct computation gives for $v\in C_c^2(0,\infty)$,
\begin{equation}\label{eq:5.3}\frac{F_3(v)}3=\int_0^\infty r(A^2+B^2)\dd
\end{equation}
and
\begin{equation}
\label{eq:5.4}\frac{E_3(v)}3=\int_0^\infty r^2(A'^2+A^2)\dd+\int_0^\infty\bigl[r^2(B'^2+B^2)+6B^2\bigr]\dd =:\overline Q_0[A] + \overline Q_2[B].\end{equation}
We have then
\begin{lemma}\label{lem:5.2}
If $A$ and $B$ satisfy \eqref{eq:5.1} and the corresponding $v$ vanishes at infinity, then
\begin{equation}\label{eq:5.7}A(r)=T_3B(r):=\frac1{\sqrt2}\Big(B(r)-3\int_r^\infty\frac{B(s)}s\,ds\Big).\end{equation}
The map $T_3$ is continuous from the $B$-channel form domain into the $A$-channel form domain.
\end{lemma}
Indeed, the formula \eqref{eq:5.7} follows from \eqref{eq:5.1}, we apply \eqref{eq:3.8} with $a=1, 2$ to control $\|A\|_{L^2((r+r^2)dr)}$ and there holds
\[A'=\frac1{\sqrt2}\Big(B'+\frac{3B}r\Big),\]
so we can control $\overline Q_0[A]$ by $\overline Q_2[B]$. The channel operators are now
\[\overline\L_0 A =-(r^2A')'+r^2A,\quad \overline\L_2 B=-(r^2B')'+(r^2+6)B.\]
Define, with Laguerre polynomials $L_m^{(\beta)}$, 
\begin{equation*}\overline \phi_n(r)=e^{-r}L_n^{(1)}(2r),\quad \overline\psi_m(r)=r^2e^{-r}L_m^{(5)}(2r).\end{equation*}
Then
\begin{equation}\label{eq:5.9}\overline \L_0\overline \phi_n=(2n+2)r\overline \phi_n,\quad \overline \L_2\overline\psi_m =(2m+6)r\overline \psi_m,\end{equation}
and 
\begin{equation}\label{eq:5.10}\int_0^\infty r\overline \phi_n\overline\phi_k\dd=\frac{n+1}4\delta_{nk},\quad\int_0^\infty r\overline\psi_m \overline\psi_\ell \dd = \overline h_m\delta_{m\ell},\end{equation}
with $\overline h_m =\frac{(m+1)(m+2)(m+3)(m+4)(m+5)}{64}$.

\begin{lemma}\label{lem:5.3}
For every $m\ge0$, there holds
\begin{equation}\label{eq:5.12}\sqrt2\,T_3\overline\psi_m = -\frac{3(m+3)}2\sum_{\substack{0\le n\le m\\n\equiv m\ (\mathrm{mod}\ 2)}}\overline \phi_n+\frac{(m+1)(m+2)}4 \overline\phi_{m+2}.\end{equation}
\end{lemma}
Here we use the generating functions
\begin{equation*}\sum_{n\ge0}L_n^{(1)}(2r)z^n=(1-z)^{-2}\exp\Big({-\frac{2rz}{1-z}}\Big),\quad \sum_{m\ge0}L_m^{(5)}(2r)z^m=(1-z)^{-6}\exp\Big({-\frac{2rz}{1-z}}\Big)\end{equation*}
and 
\[\sum_{n \ge 0} \overline\phi_n(r)z^n = (1-z)^{-2}e^{-a(z)r}, \quad \sum_{m\ge0}\overline\psi_m(r)z^m=r^2(1-z)^{-6}e^{-a(z)r},\]
where $a(z) = \frac{1+z}{1-z}$. Using the elementary integrals as above, there holds
\begin{equation*}
\sqrt2 \sum_{m\ge0} T_3\overline \psi_m(r)z^m=(1-z)^{-6}e^{-a(z)r}\Big[r^2-\frac{3r}{a(z)}-\frac3{a(z)^2}\Big].\end{equation*}
Summation of the right-hand side of \eqref{eq:5.12}, again separating parity tails, the coefficient comparison proves the Lemma.

Expand $B=\sum b_m\overline\psi_m$ and define
\begin{equation}\label{eq:5.14} s_n:=\sum_{\substack{m\ge n\\m\equiv n\ (\mathrm{mod}\ 2)}}\frac{3(m+3)}2b_m, \quad \mbox{hence }\; b_n=\frac2{3(n+3)}(s_n-s_{n+2}).\end{equation}
If $A=\sum a_n\overline \phi_n$, then
\begin{equation}\label{eq:5.15}\sqrt2\,a_n=\overline\kappa_ns_{n-2}-(1+\overline\kappa_n)s_n,\quad\overline\kappa_n=\frac{n(n-1)}{6(n+1)}.\end{equation}
\begin{proposition}\label{prop:5.4}
Finite linear combinations of $\{\overline\psi_m\}$ are a core for the $B$-channel quadratic form in $L^2(\R_+, rdr)$. Under $T_3$, the corresponding images converge in the $A$-channel form norm. Hence finite Laguerre expansions are dense in the original $(E_3+F_3)^{1/2}$ energy space.
\end{proposition}
For the Jacobi chain, again the parity sectors decouple. The odd sector has bottom quotient at least $4$, while \eqref{eq:1.2} implies $\C_3<4$. Similarly to the two dimensional case, we can claim that the ground state lies in the even chain. Set $x_j=s_{2j}$ for $j\ge0$,
and define
\begin{equation*}\overline{\mathcal E}_j(x,y)= \overline p_j x^2-2\overline q_jxy+ \overline r_jy^2,\end{equation*}
where $n=2j$,
\begin{equation}\label{eq:6.3}\overline p_j =\frac{(n+1)(n+2)(n^2+7n+14)}{96(n+3)},\quad \overline q_j = \frac{(n+1)(n+2)(n+4)(n+5)}{96(n+3)}\end{equation}
and 
\begin{equation}\label{eq:6.5}\overline r_j =\frac{(n+4)(n+5)(n^2+5n+8)}{96(n+3)}.\end{equation}
The determinant is
\begin{equation}\label{eq:6.6} \overline p_j \overline r_j - \overline q_j^2=\frac{(n+1)(n+2)(n+4)(n+5)}{1152}>0.\end{equation}
\begin{theorem}
\label{thm:6.1}
For $2<\alpha<4$, we get
\begin{equation}\label{eq:6.7}\frac{G_{3,\alpha,\rm even}(v)}3=\frac{2-\alpha}8x_0^2+\sum_{j=0}^\infty(4j+6-\alpha)\overline{\mathcal E}_j(x_j,x_{j+1}).\end{equation}
Consequently
\begin{equation}\label{eq:6.8} \C_3=\inf_{x\ne0}\frac{\frac14x_0^2+\sum_{j\ge0}(4j+6)\overline{\mathcal E_j}(x_j,x_{j+1})}{\frac18x_0^2+\sum_{j\ge0}\overline{\mathcal E_j}(x_j,x_{j+1})}.\end{equation}
\end{theorem}
\begin{proof} By \eqref{eq:5.9}--\eqref{eq:5.10}, there holds
\[\frac{G_{3,\alpha,\rm even}(v)}3=\sum_{n\ge0}(2n+2-\alpha)\frac{n+1}4a_n^2+\sum_{n\ge0}(2n+6-\alpha)\overline h_n b_n^2.\]
For the even sector, $a_0=-\frac{1}{\sqrt2}x_0$, which gives the boundary term $\frac{1}{8}(2-\alpha)x_0^2$. Substitution of \eqref{eq:5.14} and \eqref{eq:5.15} gives \eqref{eq:6.7}. Moreover, put
\[w_n=\frac{(n+1)(n+2)(n+4)(n+5)}{144(n+3)},\quad\sigma_n =\frac{(n+1)(n+2)}{6(n+3)}.\]
The pair $(b_n,a_{n+2})$ contributes
\begin{equation}\label{eq:C.1}(2n+6-\alpha)\Big\{w_n(x_j-x_{j+1})^2+\frac{n+3}8\bigl[\sigma_n x_j-(1+\sigma_n)x_{j+1}\bigr]^2\Big\}.\end{equation}
Expanding \eqref{eq:C.1} gives \eqref{eq:6.3}--\eqref{eq:6.5}, hence \eqref{eq:6.6}. Proposition~\ref{prop:5.4} justifies the form-domain limit. Separating the $\alpha=0$ and $\alpha$-dependent parts gives \eqref{eq:6.8}.
\end{proof}

Let $\overline\lambda_j =4j+6$ and $\beta_3(\alpha)=\frac{\alpha-2}8$. Then for $\alpha<\overline\lambda_0 = 6$, the tail energy
\begin{equation*}
\overline R_j(\alpha):=\inf_{\substack{x_j=1\\x\text{ in the chain form domain}}}\sum_{k=j}^\infty(\overline \lambda_k -\alpha)\overline{\mathcal E}_k(x_k,x_{k+1}) > 0.
\end{equation*}
The iterative formula \eqref{eq:7.2} remains valid if we replace $p_j, q_j, \lambda_j, R_j$ by $\overline p_j, \overline q_j, \overline\lambda_j, \overline R_j$. Therefore, for $d_j = \overline\lambda_j-\alpha =4j+6-\alpha$, we get the continued fraction expression
\begin{equation}\label{eq:7.4}
\overline R_0(\alpha) = d_0\overline p_0 -\cfrac{d_0^2\overline q_0^2}{d_0\overline r_0 + d_1\overline p_1 - \cfrac{d_1^2\overline q_1^2}{d_1\overline r_1+ d_2 \overline p_2-\cdots}}.
\end{equation}
We conclude 
\begin{theorem}\label{thm:7.2bis}
The constant $\C_3=2\H_3$ is characterized by $\overline R_0(\C_3)=\beta_3(\C_3)$. This equation has at most one solution in $(0, 6)$, hence the solution in Theorem~\ref{thm:main} is unique.
\end{theorem}
The approximation using the backward Riccati recursion, if we replace $R_J^{\pm,J}, p_J, \lambda_J$ by $\overline R_J^{\pm,J}, \overline p_J, \overline\lambda_J$, and 
$\overline \alpha_J^\pm$ be the roots of $\overline R_0^{\pm,J}(\alpha)=\beta_3(\alpha)$. So we can claim that $\overline \alpha_J^\pm$ are adjacent sequences which converge to $\C_3$. Finally, the extremizer is also unique up to a scaling with $\overline x_0 = 1$ and 
\begin{equation*}
\overline x_{j+1}=\frac{(\overline \lambda_j-\C_3)\overline q_j}{(\overline \lambda_j-\C_3)\overline r_j + \overline R_{j+1}(\C_3)}\overline x_j, \quad \forall\; j \in \N.
\end{equation*}
\appendix

\section{Finite tridiagonal implementation}
For a rational target $\alpha$, the finite matrix representing $\frac{G_{2,\alpha, \rm even}}2$
on the coordinates \((x_0,\dots,x_J)\) is tridiagonal. With $c_j=\lambda_j-\alpha=4j+5-\alpha$, the offdiagonal entries are
\[
E_j=-c_jq_j,\quad 0\le j<J.
\]
The diagonal entries are $D_0=\frac{1-\alpha}{2}+c_0p_0$, and for $1\le j<J$, $D_j=c_{j-1}r_{j-1}+c_jp_j$. For the free-tail matrix one takes $D_J=c_{J-1}r_{J-1}$, while for the Dirichlet matrix one takes $D_J=c_{J-1}r_{J-1}+c_Jp_J$. Here $p_j,q_j,r_j$ are the coefficients defined in \eqref{eq:4.3}--\eqref{eq:4.5}.

With this notation, the exact positivity test is given by the scalar $LDL^T$ recursion
\[d_0=D_0,\quad d_j=D_j-\frac{E_{j-1}^2}{d_{j-1}}.\]
This is algebraically the finite continued fraction.

Similar argument holds for dimension three, we need only to replace $\lambda_j, p_j, q_j, r_j$ by $\overline\lambda_j, \overline p_j, \overline q_j, \overline r_j$, so we omit the details.

\bigskip
\noindent
{\bf AI assistance statement.} During the preparation of this manuscript, the authors used AI-assisted tools for exploratory discussions, mathematical checking, and editing. The authors independently verified the mathematical content and take full responsibility for the work.

\end{document}